\documentclass[a4paper,12pt]{article}
\usepackage{amssymb}
\usepackage{amsmath}
\usepackage{epsfig}
\usepackage{color}

\newcommand{\ZZ}{{\mathbb Z}}
\newcommand{\G}{\Gamma}

\newcommand{\la}{\langle}
\newcommand{\ra}{\rangle}
\newcommand{\qed}{\hfill\hbox{\rule{3pt}{6pt}} \medskip}
\newcommand{\proof}{{\sc Proof. }}

\newcommand{\GP}{\hbox{GP}}
\newcommand{\Aut}{\mathrm{Aut}}

\newtheorem{theorem}{Theorem}[section]
\newtheorem{lemma}[theorem]{Lemma}
\newtheorem{corollary}[theorem]{Corollary}
\newtheorem{proposition}[theorem]{Proposition}

\newtheorem{remark}[theorem]{Remark}

\newtheorem{problem}[theorem]{Problem}

\newtheorem{conjecture}[theorem]{Conjecture}
\begin{document}
\begin{center}
{\bf\Large $\ell$-DISTANCE-BALANCED GRAPHS} \\ [+4ex]
\v{S}tefko Miklavi\v{c}{\small$^{a, b, c, }$}\footnotemark,\
Primo\v z \v Sparl{\small$^{b, c, d, }$}\footnotemark,  
\\ [+2ex]
{\it \small 
$^a$University of Primorska, FAMNIT, Glagolja\v ska 8, 6000 Koper, Slovenia\\
$^b$University of Primorska, IAM, Muzejski trg 2, 6000 Koper, Slovenia\\
$^c$IMFM, Jadranska 19, 1000 Ljubljana, Slovenia\\
$^d$University of Ljubljana, Faculty of Education, Kardeljeva plo\v s\v cad 16, 1000 Ljubljana, Slovenia}
\end{center}

\addtocounter{footnote}{-1}
\footnotetext{The author acknowledges the financial support from the Slovenian Research Agency (research core funding No. P1-0285 and research projects N1-0032, N1-0038, J1-6720, J1-7051).}
\addtocounter{footnote}{1} \footnotetext{
The author acknowledges the financial support from the Slovenian Research Agency (research core funding No. P1-0285 and research projects N1-0038, J1-6720, J1-7051).\\

Email addresses: 
stefko.miklavic@upr.si (\v Stefko Miklavi\v c),
primoz.sparl@pef.uni-lj.si (Primo\v z \v Sparl).
}

%


\begin{abstract}
Let $\ell$ denote a positive integer.
A connected graph $\G$ of diameter at least $\ell$ is said to be $\ell${\it -distance-balanced} whenever for any pair of vertices $u,v$ of
$\G$ such that $d(u,v)=\ell$, the number of vertices closer to $u$ than to $v$ is equal to the number of
vertices closer to $v$ than to $u$. 
In this paper we present some basic properties of $\ell$-distance-balanced graphs and study in more detail $\ell$-distance-balanced graphs of diameter at most $3$. We also investigate the $\ell$-distance-balanced property of some well known families of graphs such as the generalized Petersen graphs.
\end{abstract}

\section{Introduction}
\label{sec:intro}

Throughout this paper, all graphs are simple (without loops and multiple edges), undirected, finite and connected.
Given a graph $\G$ let $V(\G)$ and $E(\G)$ denote its vertex set and edge set, respectively.
For $u,v \in V(\G)$ we denote the distance between $u$ and $v$ in $\Gamma$ by $d_\Gamma(u,v)$ (or simply $d(u,v)$ if the graph $\Gamma$ is clear from the context).
Furthermore, for any nonnegative integer $i$ and $u \in V(\G)$ let $N_i(u)=\{v \in V(\G) \mid d(u,v)=i\}$
(we abbreviate $N(u) = N_1(u)$).
For $S \subseteq V(\G)$ the subgraph of $\G$ induced by $S$ is denoted by $\la S \ra$
(we abbreviate $\G-S = \la V(\G) \setminus S \ra$). 

For any pair of vertices $u,v \in V(\G)$ we let $W_{uv}$ be the set of vertices of $\G$ that are closer to $u$ than to $v$, that is
$$
  W_{uv} = \{w \in V(\G) \mid d(u,w) < d(v,w) \}.
$$
The pair $u, v$ is said to be {\em balanced} if $|W_{uv}| = |W_{vu}|$ and is {\em non-balanced} otherwise. Let $\ell$ denote a positive integer.
A connected graph $\G$ of diameter at least $\ell$ is said to be $\ell$-{\em distance-balanced} whenever any pair of vertices $u,v \in V(\G)$ at distance $\ell$ is balanced, that is, if for any $u, v \in V(\G)$ such that $d(u,v)=\ell$ we have
$$
 |W_{uv}| = |W_{vu}|.
$$
A connected graph $\G$ is said to be {\em highly distance-balanced} if it is $\ell$-distance-balanced for every $1 \le \ell \le D$, where $D$ is the diameter of $\Gamma$.  

The $\ell$-distance-balanced graphs are a natural generalization of the so-called {\em distance-balanced} graphs~\cite{JKR}. They were first defined by Bo\v{s}tjan Frelih in his PhD disertation \cite{F}, where  $2$-distance-balanced graphs were studied in more detail. In particular, $2$-distance-balanced graphs which are not 2-connected were characterized, and $2$-distance-balanced graphs were studied with respect to various graph products. 
The $\ell$-distance-balanced graphs were also the main topic of the paper~\cite{FPK}. However, some of the stated results do not hold while some are given without proof. We comment on two of these problems later (see Remarks~\ref{remark1} and~\ref{remark2}).

On the other hand,  
distance-balanced graphs have been extensively studied, see \cite{BCPSSS, Ha, HN84, IKM, JKR, KMMM06, KMMM2, KM, MS}. We also point out that every distance-regular graph~\cite{BCN} is highly distance-balanced. The opposite is of course not true. For instance, the generalized Petersen graph $GP(7,2)$ (see the definition at the end of this section) is highly distance-balanced (see Table~\ref{tab:GP}) but it can be easily seen that it is not distance-regular.

Let $G$ be a group and let $S \subset G$ be an inverse closed subset (that is $S = S^{-1}$) not containing the identity. Then the {\em Cayley graph} $\mathrm{Cay}(G;S)$ is defined to be the graph with vertex set $G$ in which $g \in G$ is adjacent to $h \in G$ whenever $g^{-1}h \in S$. 

Let $n \ge 3$ be a positive integer, and let 
$1 \leq k < n/2$.
The generalized Petersen graph $\GP(n, k)$ 
is defined to have the following vertex set and edge set:
\begin{eqnarray}
\label{genpet}
  V(\GP(n,k)) &=& \{ u_i \mid i \in \ZZ_n \}\cup \{v_i \mid i \in \ZZ_n \},  \nonumber \\
  E(\GP(n, k)) &=& \{u_iu_{i+1} \mid i \in \ZZ_n\}\cup
    \{ v_{i}v_{i+k}  \mid i \in \ZZ_n\}\cup
    \{  u_iv_i \mid i \in \ZZ_n\}.
\end{eqnarray}
\noindent
The edges of the form $u_iu_{i+1}$ are called {\em outer edges}, edges of the form $v_iv_{i+k}$ are called {\em inner edges},
and edges of the form $u_iv_i$ are called {\em spokes}.
Note that $\GP(n, k)$ is cubic, and
that it is bipartite precisely when $n$ is even and $k$ is odd.
It is easy to see that $\GP(n,k)\cong\GP(n,n-k)$. Furthermore, if
the multiplicative inverse $k^{-1}$ of $k$ exists in $\ZZ_n$, then
$\GP(n,k) \cong \GP(n, k^{-1})$.

In this paper we first study basic properties of $\ell$-distance-balanced graphs.
We also give examples of these graphs.
In Section \ref{sec:D<4} we study $\ell$-distance-balanced graphs with diameter at most $3$. 
In Section \ref{sec:GP} we study the $\ell$-distance-balanced property of the generalized Petersen graphs.

\section{Basic properties and examples}
\label{sec:basic}

In this section we present some basic properties of $\ell$-distance-balanced graphs and give various examples of such graphs. We first state a fairly straightforward but useful observation and its corollary.

\begin{lemma}
\label{le:swap}
Let $\Gamma$ be a connected graph and $u,v \in V(\Gamma)$. If some $\alpha \in \Aut(\Gamma)$ interchanges $u$ and $v$, then the pair $u,v$ is balanced. 
\end{lemma}
\proof
This follows from the fact that for any automorphism $\alpha \in \Aut(\Gamma)$ we have $\alpha(W_{uv}) = W_{\alpha(u) \alpha(v)}$, and so in the case that $\alpha$ interchanges $u$ and $v$ we obtain $\alpha(W_{uv}) = W_{vu}$. As $\alpha$ is a bijection we thus get $|W_{uv}| = |W_{vu}|$. 
\qed

\begin{corollary}
\label{cor:swap}
Let $\G$ be a connected graph such that for each pair of vertices $u, v \in V(\G)$ there exists an automorphism of $\G$ interchanging $u$ and $v$. Then $\G$ is highly distance-balanced.
\end{corollary}
\begin{remark}
\label{remark1}
A similar result as in Corollary~\ref{cor:swap} is given in~\cite[Proposition 2.10]{FPK} but without a proof. It is then stated (in \cite[Corollary 2.12]{FPK}) that this forces every graph $\G$, in which for every pair of vertices at distance $\ell$ there is an automorphism of $\G$ mapping one to the other, to be $\ell$-distance-balanced. This however does not hold, since it would imply that every vertex-transitive graph is highly distance-balanced. But, for example, the generalized Petersen graph $GP(16,7)$ with diameter $5$ is vertex-transitive (as $7^2$ is congruent to $1$ modulo $16$), but it is not $4$-distance-balanced (see Table~\ref{tab:GP}).
\end{remark}

Corollary \ref{cor:swap} for instance implies that the Petersen graph $GP(5,2)$ is highly distance-balanced since it is distance-transitive. Moreover, the corollary gives rise to infinitely many highly distance-balanced graphs. 
In particular, every Cayley graph of an abelian group is highly distance-balanced.

\begin{proposition}
\label{pro:abelian}
Let $A$ be a finite abelian group and let $S \subset A$ be an inverse closed subset of $A$ not containing the identity $1$. If $\langle S \rangle  = A$ then the Cayley graph $\mathrm{Cay}(A;S)$ is highly distance-balanced. 
\end{proposition}
\proof
Since the graph $\Gamma = \mathrm{Cay}(A;S)$ is vertex transitive it suffices to prove that there exists no $a \in A$ such that the pair 
$1$, $a$ is non-balanced. Observe that, since $A$ is abelian, the permutation $\tau$ of $A$, mapping each element to its inverse, is an automorphism of $\Gamma$. 
Namely, for any pair $a$, $as$ of adjacent vertices of $\Gamma$ their images $a^{-1}$ and $s^{-1}a^{-1} = a^{-1}s^{-1}$ are adjacent as $S = S^{-1}$. 
Likewise, for any $a \in A$ the permutation $t_a$, mapping each $b \in A$ to $ab$, is clearly an automorphism of $\Gamma$. Since the product $t_a \tau$ interchanges the vertices $1$ and $a$, Corollary~\ref{cor:swap} implies that $\Gamma$ is highly distance-balanced. 
\qed

The above proposition implies that cycles, being Cayley graphs of cyclic groups, are highly distance-balanced. One of the next natural families of graphs that could
possibly provide interesting examples and nonexamples of $\ell$-distance-balanced graphs is the family of cubic graphs. Within this family the well known family of generalized Petersen graphs might be a good place to start the investigation. As we will see in Section~\ref{sec:GP} the problem of determining all $\ell$ such that, for given $n \geq 3$ and $1 \leq k < n/2$, the generalized Petersen graph $GP(n,k)$ is $\ell$-distance-balanced, does not seem to be easy. Of course, for some pairs of $n$ and $k$, the problem is very easy. For instance, the prisms, being Cayley graphs of abelian groups, are not very interesting.

\begin{corollary}
\label{cor:prism}
Let $n \geq 3$ be an integer. Then the prism $GP(n,1)$ is highly distance-balanced.
\end{corollary}
\proof
Since the genarlized Petersen graph $GP(n,1)$ is isomorphic to the Cayley graph $\mathrm{Cay}(\ZZ_n \times \ZZ_2 ; \{(1,0), (-1,0), (0,1)\})$, it is 
highly distance-balanced by Proposition~\ref{pro:abelian}. \qed

\begin{figure}[h]
\begin{center}
\includegraphics[scale=0.7]{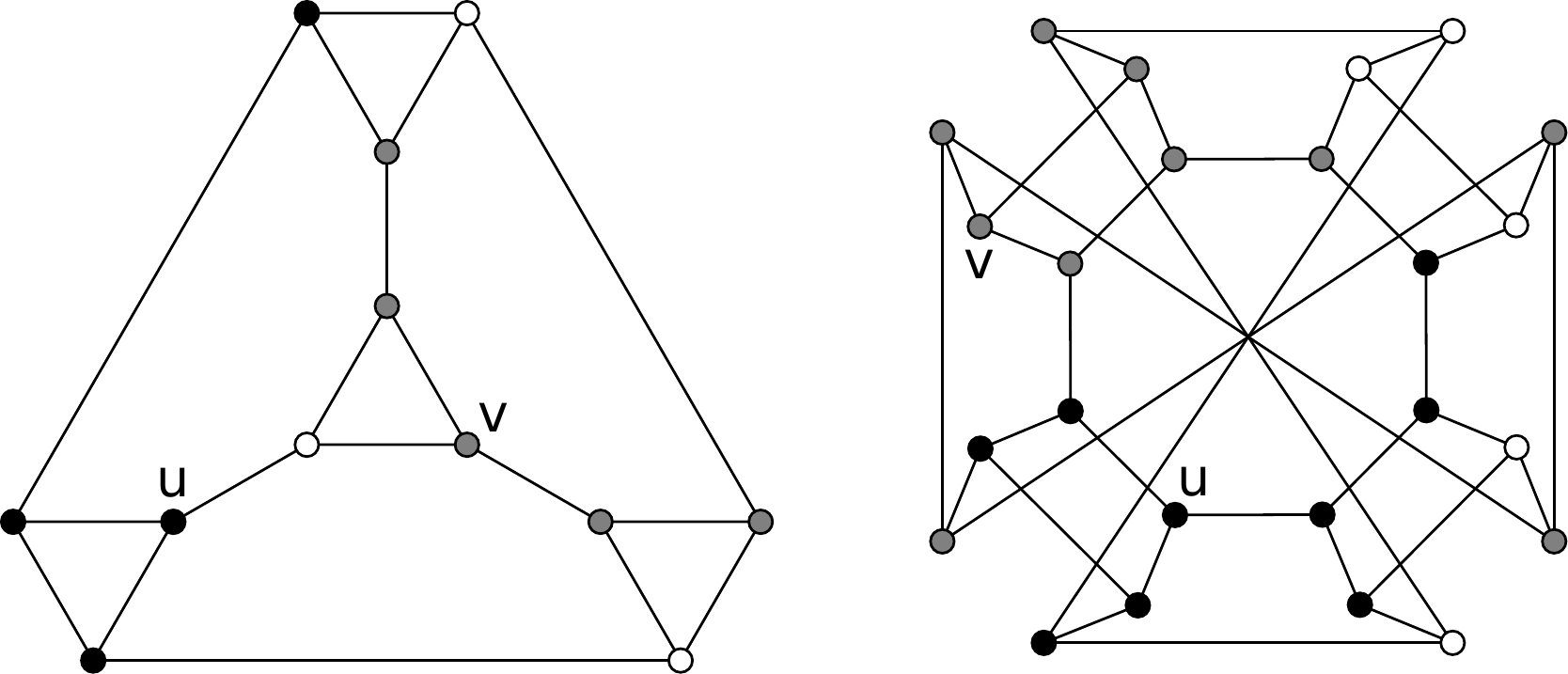}
\end{center}
\caption{The Cayley graph $\mathrm{Cay}(A_4; \{(1\,2\,3), (1\,3\,2), (1\,2)(3\,4)\})$ is not $2$-distance-balanced and the Cayley graph $\mathrm{Cay}(S_4 ; \{(1\,2), (2\,4), (1\,2)(3\,4)\})$ is not $3$-distance-balanced.}
\label{fig:examples}
\end{figure}

As we will see in Section \ref{sec:GP}, the generalized Petersen graphs $GP(n,k)$ with $k \geq 2$ are much more interesting. However, before we turn 
our attention to these graphs let us mention a few more interesting examples. 
The Cayley graph $\mathrm{Cay}(A_4; S)$, where $S=\{(1\,2\,3), (1\,3\,2), (1\,2)(3\,4)\})$ (the truncation of the tetrahedron), which is a cubic graph of diameter $3$, 
is $1$-distance-balanced and $3$-distance-balanced, but it is not $2$-distance-balanced. Namely, the pair $u = \mathrm{id}$ and $v = (1\,3\,4)$ is not balanced as we get $|W_{uv}| = 4$ and $|W_{vu}| = 5$ (see Figure~\ref{fig:examples}). Similarly, the Cayley graph $\mathrm{Cay}(S_4 ; \{(1\,2), (2\,4), (1\,2)(3\,4)\})$, which is a cubic graph of diameter $4$, is $1$-distance-balanced and $2$-distance-balanced but is not $3$-distance-balanced (see Figure~\ref{fig:examples}) nor $4$-distance-balanced since none of the pairs $\mathrm{id}$, $(1\,4\,3)$, nor $\mathrm{id}$, $(1\,2\,4\,3)$ is balanced.
It thus seems that already with cubic graphs the situation regarding $\ell$-distance-balancedness is quite interesting.

\smallskip
Another interesting family of graphs are distance degree regular graphs (first introduced in~\cite{HN84} and later called strongly distance-balanced graphs in~\cite{KMMM06}). 
A connected graph $\G$ with diameter $D$ is called {\em distance degree regular}, whenever 
$|N_i(u)|= |N_i(v)|$ for any two vertices $u,v$ of $\G$ and for any $0 \le i \le D$. It was shown in~\cite{KMMM06} that every distance degree regular graph is automatically distance-balanced. As we now show, it is also $2$-distance-balanced provided it is bipartite.

\begin{theorem}
\label{ddr_thm1}
Let $\G$ denote a bipartite distance degree regular graph. Then $\G$ is $1$- and $2$-distance-balanced.
\end{theorem}
\proof
By the above remark we only need to show that $\Gamma$ is $2$-distance-balanced. Pick vertices $u,v$ of $\G$ such that $d(u,v)=2$. We show that then $|W_{uv}| = |W_{vu}|$. 
Observe first that since $\G$ is bipartite we have 
\begin{equation}
\label{ddr_eq1}
  W_{uv} = \bigcup_{i=1}^{D-1} (N_{i-1}(u) \cap N_{i+1}(v)), \qquad
  W_{vu} = \bigcup_{i=1}^{D-1} (N_{i-1}(v) \cap N_{i+1}(u)).
\end{equation}
To prove the theorem it thus suffices to verify that for each $1 \leq i \leq D-1$ the equality 
$|N_{i-1}(u) \cap N_{i+1}(v)|=|N_{i-1}(v) \cap N_{i+1}(u)|$ holds.
We show this using induction on $i$. Obviously, $|N_0(u) \cap N_2(v)|=|N_0(v) \cap N_2(u)|=1$. Since $d(u,v) = 2$ and $\Gamma$ is bipartite we have that $N_1(u) = (N_1(u) \cap N_1(v)) \cup (N_1(u) \cap N_3(v))$ and $N_1(v) = (N_1(u) \cap N_1(v)) \cup (N_1(v) \cap N_3(u))$. Thus, since $|N_1(u)| = |N_1(v)|$ we get $|N_1(u) \cap N_3(v)|=|N_1(v) \cap N_3(u)|$.

Suppose now that for some $2 \le k \le D-2$ we have that 
$|N_{j-1}(u) \cap N_{j+1}(v)|=|N_{j-1}(v) \cap N_{j+1}(u)|$ for
each $1 \le j \le k$. Observe that $N_k(u)$ is a disjoint union of 
$N_k(u) \cap N_{k+2}(v)$, $N_k(u) \cap N_k(v)$ and $N_k(u) \cap N_{k-2}(v)$.
Similarly, $N_k(v)$ is a disjoint union of 
$N_k(v) \cap N_{k+2}(u)$, $N_k(v) \cap N_k(u)$ and $N_k(v) \cap N_{k-2}(u)$.
Since $|N_k(u)| = |N_k(v)|$ we thus get
$$
  |N_k(u) \cap N_{k+2}(v)| + |N_k(u) \cap N_{k-2}(v)| = 
  |N_k(v) \cap N_{k+2}(u)| + |N_k(v) \cap N_{k-2}(u)|.
$$  
Since, by induction hypothesis, $|N_k(u) \cap N_{k-2}(v)| = |N_k(v) \cap N_{k-2}(u)|$ holds, we thus obtain 
$|N_k(u) \cap N_{k+2}(v)| = |N_k(v) \cap N_{k+2}(u)|$, which completes the induction step.
\qed

\begin{corollary}
\label{ddr_cor2}
Let $\G$ denote a connected bipartite vertex transitive graph. Then $\G$ is $2$-distance-balanced.
In particular, every bipartite connected Cayley graph is $2$-distance-balanced.
\end{corollary}
\proof
Observe that every vertex transitive graph is clearly distance degree regular. 
The result now follows immediately from Theorem \ref{ddr_thm1}. \qed

We remark that the result of Theorem~\ref{ddr_thm1} (as well as Corollary~\ref{ddr_cor2}) cannot be extended to nonbipartite distance degree regular graphs. Namely, the truncation of the tetrahedron is vertex-transitive (being a Cayley graph) and as such is distance degree regular but is not $2$-distance-balanced as was indicated on Figure~\ref{fig:examples}.

\section{Graphs of diameter at most $3$}
\label{sec:D<4}

In this section we study graphs with diameter $2$ or $3$ (the graphs of diameter $1$, i. e. the complete graphs, are of course $1$-distance-balanced). Assume first that $\G$ has diameter $2$.
By \cite[Corollary 2.3]{JKR}, $\G$ is $1$-distance-balanced if and only if it is regular. It thus remains to determine when $\G$ is $2$-distance-balanced. To do so we first need some more terminology. 

Let $\G_1$ and $\G_2$ be graphs with disjoint vertex sets $V_1$ and $V_2$ 
and edge sets $E_1$ and $E_2$. The {\em union} $\G=\G_1 \cup \G_2$ of graphs $\G_1$ and $\G_2$ is the graph with vertex set 
$V=V_1 \cup V_2$ and edge set $E=E_1 \cup E_2$.
The {\em join} $\G=\G_1+\G_2$ of graphs $\G_1$ and $\G_2$ is the graph $\G_1 \cup \G_2$ together with all the edges joining 
$V_1$ and $V_2$. Note that the join operation is both commutative and associative. We can thus speak of the graph $\G=\G_1 + \G_2 + \cdots + \G_t$ whenever the graphs $\G_i$ have pairwise disjoint vertex- and edge-sets. We call the graphs $\G_i \; (1 \le i \le t)$ the {\em components of} the join $\G$.
We remark that if a graph $\G$ is a join of at least two graphs then clearly the diameter of $\G$ is at most $2$.

\begin{theorem}
\label{the:D=2_2db}
Let $\G$ be a graph with diameter $2$. Then $\G$ is $2$-distance-balanced if and only if it is a join of regular graphs, that is $\G = \G_1 + \G_2 + \cdots + \G_t$, where $t \geq 1$ and each of $\G_i$ is a regular graph.
\end{theorem}
\proof
Suppose first that $\G = \G_1 + \G_2 + \cdots + \G_t$, where each $\G_i$, $1 \leq i \leq t$ is a regular graph. By definition of a join of graphs any two vertices from different components $\G_i$ are adjacent. Thus, if $u, v$ are any two vertices of $\G$ at distance $2$ then there exists some $1 \leq i \leq t$ such that $u$ and $v$ both belong to $\G_i$. Since both $u$ and $v$ are adjacent to all the vertices that are not in $\G_i$ and $\G$ is of diameter $2$, it is clear that $W_{uv}$ consists of $u$ and the neighbours of $u$ in $\G_i$, which are not neighbours of $v$. Similarly, $W_{vu}$ consists of $v$ and the neighbours of $v$ in $\G_i$, which are not neighbours of $u$. Since $\G_i$ is regular, this implies that $\G$ is $2$-distance-balanced.

Suppose now that $\G$ is $2$-distance-balanced and take an arbitrary pair of vertices $u, v$ of $\G$ at distance $2$. Let $C=N(u) \cap N(v)$. Since $\G$ is of diameter $2$, we get
$$
  W_{uv} = \{u\} \cup (N(u) \setminus C), \qquad W_{vu} = \{v\} \cup (N(v) \setminus C),
$$
and so $|W_{uv}| = |W_{vu}|$ implies that $u$ and $v$ have the same valence. Therefore, any two vertices at distance $2$ in $\G$ have the same valence. Let now $k_1 < k_2 < \cdots < k_t$ be all possible degrees of vertices of $\G$. If $t = 1$ the graph $\G$ is regular, so the proof is complete. Suppose then that $t \ge 2$, let $V_i = \{w \in V(\G) \mid |N(w)|=k_i\}$ for all $1 \le i \le t$ and let $\G_i$ denote the subgraph of $\G$, induced on $V_i$. To complete the proof we need to show that any two vertices from different sets $V_i$ are adjacent and that each $\G_i$ is a regular graph. That the former is true follows from the fact that for $u \in V_i$ and $v \in V_j$ with $1 \le i < j \le t$ the valencies of $u$ and $v$ are different, and so they cannot be at distance $2$ by the above argument (recall that $\G$ has diameter $2$). That any two vertices from the same set $V_i$ have the same valence within $\G_i$ is now clear since they have the same valence within $\G$ and they are both adjacent to all the vertices $w \in V(\G) \setminus V_i$.\qed

The following theorem is an immediate corollary of the above theorem and \cite[Corollary 2.3]{JKR}.

\begin{theorem}
\label{the:D=2}
Let $\G$ be a graph with diameter $2$. Then $\G$ is highly distance-balanced if and only if it is regular. Moreover, it is $2$-distance-balanced but not $1$-distance-balanced if and only if it is a nonregular join of at least two regular graphs. 
\end{theorem}

Let us point out the following interesting consequence of the above results. For graphs of diameter $2$ the fact that the graph in question is $1$-distance-balanced implies it is highly distance-balanced. To see that for graphs of larger diameter this does not hold in general it suffices to look at the examples from Figure~\ref{fig:examples}. There we found a $1$-distance-balanced graph of diameter $3$ which is not $2$-distance-balanced and a $1$- and $2$-distance-balanced graph of diameter $4$ that is not $3$-distance-balanced. One might think that perhaps a $1$-distance-balanced graph $\G$ is always $D$-distance-balanced where $D$ is the diameter of $\G$. However, since the graph of diameter $4$ from Figure~\ref{fig:examples} is not $4$-distance-balanced this also does not hold.
\bigskip
\begin{remark}
\label{remark2}
A similar result about $2$-distance-balanced graphs of diameter $2$ as in Theorem~\ref{the:D=2} was stated in \cite[Corollary 2.4]{FPK}. The authors do not provide a proof and claim it is a corollary of \cite[Proposition 2.2]{FPK}, which is supposed to give a necessary and sufficient condition for a graph to be $\ell$-distance-balanced. However, the condition is neither necessary nor sufficient. For instance, the generalized Petersen graph $GP(13,3)$ is of diameter $5$ and is not $3$-distance-balanced but is $4$-distance-balanced (see Table~\ref{tab:GP}). However, one can easily check that \cite[Proposition 2.2]{FPK} claims it is $3$-distance-balanced but not $4$-distance-balanced.
\end{remark}

We now turn our attention to graphs of diameter $3$. It seems that in this case the general situation is too complicated, so we restrict our consideration to bipartite graphs. Recall that a graph having vertices of two different degrees is called {\em biregular}.

\begin{proposition}
\label{D=3_pro1}
Let $\G$ be a bipartite graph of diameter $3$ and bipartition sets $X$, $Y$. Then $\G$ is $1$-distance-balanced 
if and only if it is regular, or it is biregular with $|N(u)| = |Y|/2$ and $|N(v)| = |X|/2$ for every $u \in X$ and $v \in Y$.
\end{proposition}
\proof
Observe that, since $\G$ is bipartite and of diameter $3$, any two vertices of $X$ (or of $Y$) are at distance $2$, and so for all $u \in X$ and $v \in Y$ we have
\begin{equation}
\label{eq:D=3thm1}
	\begin{array}{c}
		N_2(u) = X \setminus \{u\},\quad N_2(v) = Y \setminus \{v\},\\ 
		X = N(v) \cup N_3(v)\quad \text{and}\quad Y = N(u) \cup N_3(u).
      \end{array}
\end{equation}
Let now $u \in X$ and $v \in Y$ be adjacent vertices of $\G$. Since $\G$ is bipartite of diameter $3$, the sets $W_{uv}$ and $W_{vu}$ are given by 
\begin{equation}
\label{D=3_eq3}
\begin{split}
  W_{uv} &= \{u\} \cup (N(u) \cap N_2(v)) \cup (N_2(u) \cap N_3(v)), \\ 
  W_{vu} &= \{v\} \cup (N(v) \cap N_2(u)) \cup (N_2(v) \cap N_3(u)).
\end{split}
\end{equation}
By (\ref{eq:D=3thm1}) we have that $N(u) \cap N_2(v) = N(u) \setminus \{v\}$ and $N(v) \cap N_2(u) = N(v) \setminus \{u\}$. Moreover, $N_2(u) \cap N_3(v) = X \setminus N(v)$ and $N_2(v) \cap N_3(u) = Y \setminus N(u)$. It is thus clear that the pair $u,v$ is balanced if and only if $|N(u)|-1 + |X| - |N(v)| = |N(v)|-1 + |Y| - |N(u)|$ that is 
\begin{equation}\label{eq:XY}
	2|N(u)| + |X| = 2|N(v)| + |Y|.
\end{equation}
Observe that this equality is equivalent both to 
\begin{equation}
\label{eq:XY2}
	|N(v)| = |N(u)| + {|X|-|Y| \over 2},\quad \text{and}\quad |N(u)| = |N(v)| + {|Y|-|X| \over 2}.
\end{equation}

We are now ready to finally prove the proposition. Suppose first that $\G$ is $1$-distance-balanced. Then for any pair of adjacent vertices $u \in X$ and $v \in Y$ the equalities (\ref{eq:XY2}) hold, and so any two vertices of $X$, sharing a common neighbor (in $Y$), have the same degree and likewise any two vertices of $Y$, sharing a common neighbor (in $X$), have the same degree. As $\G$ is connected, all vertices of $X$ have the same degree, say $k_X$, and all vertices of $Y$ have the same degree, say $k_Y$. Counting the edges between $X$ and $Y$ in two different ways we obtain the equality 
$$
  k_X |X| = k_Y |Y|.
$$
Pluging this into one of the equalities from~(\ref{eq:XY2}) and multiplying by $|X|$ (or $|Y|$) we get $k_X = |Y|/2$ and $k_Y = |X|/2$, as claimed (note that $\G$ is regular precisely when $|X| = |Y|$).

To prove the converse suppose first that $\G$ is regular. Then $|X| = |Y|$, and so (\ref{eq:XY2}) implies that a pair of adjacent vertices $u \in X$ and $v \in V$ is balanced if and only if they have the same degree which clearly holds since $\G$ is regular. Suppose finally that $\G$ is biregular with $|N(u)|=|Y|/2$ and $|N(v)|=|X|/2$ for any $u \in X$ and $v \in Y$. It is now clear that the equalities (\ref{eq:XY2}) hold for any pair of adjacent vertices $u \in X$ and $v \in Y$, and so every such pair of vertices is balanced. This shows that $\G$ is $1$-distance-balanced. \qed

\begin{proposition}
\label{D=3_pro2}
Let $\G$ be a bipartite graph of diameter $3$. Then $\G$ is $2$-distance-balanced if and only if the vertices from the same bipartition set have the same degree.
\end{proposition}
\proof
Observe first that since $\G$ is bipartite of diameter $3$ a pair of distinct vertices $u, v$ of $\G$ is at distance $2$ if and only if they both belong to the same bipartition set. Moreover, for any such pair of vertices we have 
$$
\begin{array}{c}
  W_{uv} = \{u\} \cup (N(u) \cap N_3(v)), \  W_{vu} = \{v\} \cup (N(v) \cap N_3(u))\quad \text{and}\\
N(u) = (N(u) \cap N(v)) \cup (N(u) \cap N_3(v)), \ N(v) = (N(v) \cap N(u)) \cup (N(v) \cap N_3(u)).
\end{array}
$$
Thus the pair $u, v$ is balanced if and only if 
$$
	|N(u)| - |N(u) \cap N(v)| = |N(v)| - |N(v) \cap N(u)|,
$$
which is equivalent to $u$ and $v$ being of the same degree. The graph $\G$ is thus $2$-distance-balanced if and only if any two vertices from the same bipartition set have the same degree, which completes the proof.
\qed

\begin{proposition}
\label{D=3_pro4}
Let $\G$ be a bipartite graph of diameter $3$ with bipartition sets $X, Y$. Then $\G$ is $3$-distance-balanced 
if and only if for any pair of vertices $u \in X$ and $v \in Y$ at distance $3$ we have $2|N(u)| + |X| = 2|N(v)| + |Y|$.
\end{proposition}
\proof
Let $u,v$ be vertices of $\G$ with $d(u,v)=3$. Then $u$ and $v$ belong to different bipartition sets, and so we may assume $u \in X$ and $v \in Y$. Recall that (\ref{eq:D=3thm1}) holds, and so 
$$
  W_{uv} = \{u\} \cup N(u) \cup (X \setminus (\{u\} \cup N(v))) \ \text{and}\ W_{vu} = \{v\} \cup N(v) \cup (Y \setminus (\{v\} \cup N(u))).
$$
The result follows. \qed

Combining the above three results we obtain the following corollary and theorem.

\begin{corollary}
\label{cor:D=3}
Let $\G$ be a bipartite graph of diameter $3$ with the bipartition sets $X$ and $Y$. Then $\G$ is highly distance-balanced if and only if it is either regular or biregular with each $u \in X$ and $v \in Y$ being of degree $|Y|/2$ and $|X|/2$, respectively.
\end{corollary}

\begin{theorem}
Let $\G$ be a bipartite graph of diameter $3$ with bipartition sets $X$ and $Y$. Then precisely one of the following holds:
\begin{itemize}\itemsep = 0pt
\item[(i)] $\Gamma$ is highly distance-balanced.
\item[(ii)] $\Gamma$ is $2$-distance-balanced but not $1$-distance-balanced nor $3$-distance-balanced.
\item[(iii)] $\Gamma$ is $3$-distance-balanced but not $1$-distance-balanced nor $2$-distance-balanced.
\item[(iv)] $\Gamma$ is not $\ell$-distance-balanced  for any $1 \leq \ell \leq 3$.
\end{itemize}
\end{theorem}
\proof
Suppose first that $\G$ is $1$-distance-balanced. Then Propositions~\ref{D=3_pro1},~\ref{D=3_pro2} and~\ref{D=3_pro4} imply that $\G$ is both $2$-distance-balanced and $3$-distance-balanced.

Suppose next that $\G$ is not $1$-distance-balanced but is $2$-distance-balanced and $3$-distance-balanced. By Propositions~\ref{D=3_pro1} and~\ref{D=3_pro2} the graph $\G$ is not regular, but is biregular with the degree of each $u \in X$ being $k_X$ and the degree of each $v \in Y$ being $k_Y$. Consequently 
\begin{equation}
\label{eq:D=3theo}
	|X|k_X = |Y|k_Y,\quad \text{and so}\quad k_X = k_Y \frac{|Y|}{|X|}.
\end{equation} 
Take now $u \in X$ and $v \in Y$ such that $d(u,v) = 3$. By Proposition~\ref{D=3_pro4} we have $2k_X + |X| = 2k_Y + |Y|$. Therefore (\ref{eq:D=3theo}) implies
$$
	2k_Y(|Y| - |X|) = |X|(|Y| - |X|).
$$
Since $\G$ is not regular, $|X| \neq |Y|$, and so $k_Y = |X|/2$. Similarly we obtain $k_X = |Y|/2$. By Proposition~\ref{D=3_pro1} the graph $\G$ is $1$-distance-balanced, a contradiction. \qed

We remark that each of the four possibilities from the above theorem can indeed occur. Every regular bipartite graph of diameter $3$ (for instance, the cube graph) is highly distance-balanced, proving that item (i) is possible. Take any bipartite graph with bipartition sets $X$ and $Y$ of cardinalities $6$ and $9$, respectively, and where each $u \in X$ has degree $6$ and each $v \in Y$ has degree $4$. It is easy to see that such graphs exist, have diameter $3$ and are $2$-distance-balanced but not $3$-distance-balanced, proving that item (ii) is possible. Finally, the path of length $3$ is clearly a $3$-distance-balanced bipartite graph of diameter $3$ which is neither $2$-distance-balanced nor $1$-distance-balanced, and so item (iii) is also possible.

\begin{figure}[h]
\begin{center}
\includegraphics[scale=0.6]{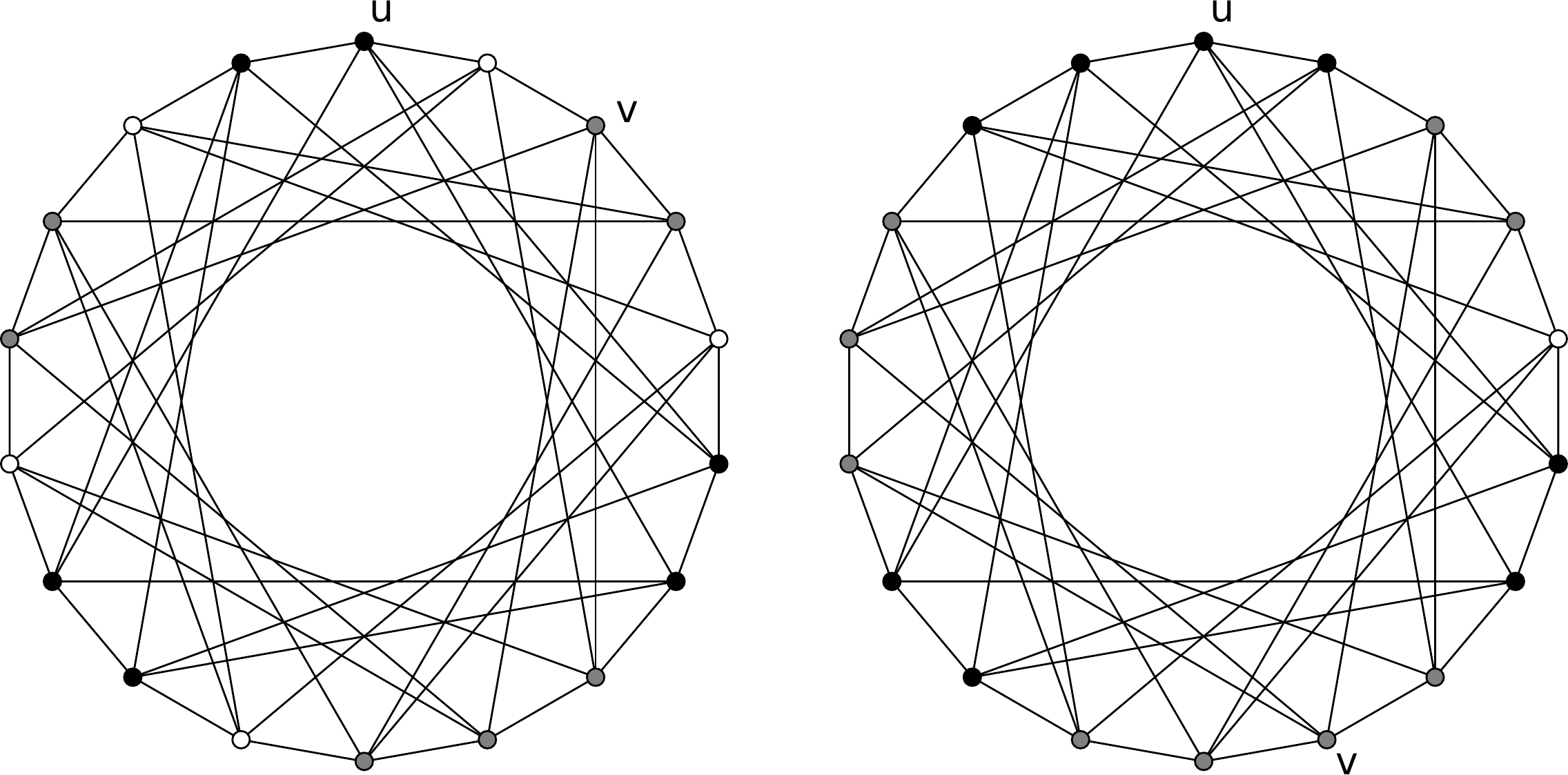}
\end{center}
\caption{The Cayley graph $\mathrm{Cay}(D_9 ; \{t, tr^2, tr^3, r^3, r^6\})$ is of diameter $3$ but is not $2$-distance-balanced nor $3$-distance-balanced.}
\label{fig:ex18}
\end{figure}

By Corollary~\ref{cor:D=3} a regular bipartitie graph of diameter $3$ is highly distance-balanced. However, if we drop the condition on bipartiteness the result no longer holds. For instance, the truncation of the tetrahedron from Figure~\ref{fig:examples} which is of course regular (being a Cayley graph) is of diameter $3$ but is not $3$-distance-balanced. This graph is $2$-distance-balanced though. However, also this need not be the case in general. For instance, the Cayley graph $\mathrm{Cay}(D_9 ; \{t, tr^2, tr^3, r^3, r^6\})$, where $D_9 = \langle t, r \mid t^2, r^9, (tr)^2 \rangle$ is of diameter $3$ but is not $2$-distance-balanced (the pair $1$, $r$ is not balanced) nor $3$-distance-balanced (the pair $1$, $r^4$ is not balanced), as can be seen on Figure~\ref{fig:ex18}.

\section{The $\ell$-distance-balanced property of generalized Petersen graphs}
\label{sec:GP}

As mentioned in Section~\ref{sec:basic} the problem of determining all $\ell$ such that a given graph is $\ell$-distance-balanced does not seem to be easy even for cubic graphs. To indicate that this might be true we investigate the well known generalized Petersen graphs and their $\ell$-distance-balancedness in this section. 

We first make the following easy but useful observation.

\begin{corollary}
\label{cor:GP}
Let $n \geq 3$ and $1 \leq k < n/2$ be integers. If the generalized Petersen graph $\Gamma = GP(n,k)$ is not $\ell$-distance-balanced for some $1 \leq \ell \leq D$, where $D$ is the diameter of $\Gamma$, then there exists $j \in \ZZ_n$ such that $d(u_0, v_j) = \ell$ and there is no automorphism of $\Gamma$ interchanging $u_0$ and $v_j$.
\end{corollary}
\proof
Observe first that the permutations $\rho = (u_0, u_1, \ldots , u_{n-1})(v_0, v_1, \ldots , v_{n-1})$ and $\tau$, where $\tau(u_i) = u_{-i}$ and $\tau(v_i) = v_{-i}$ for all $i \in \ZZ_n$, are automorphisms of $\Gamma$. By Lemma~\ref{le:swap} it thus follows that each pair $u_i$, $u_j$ and each pair $v_i$, $v_j$ is balanced. The result now follows immediately from Lemma~\ref{le:swap}.\qed

Using Corollary~\ref{cor:GP} and a suitable software package such as {\sc Magma}~\cite{Mag} one may now easily compute all the values $\ell$ for which a given $GP(n,k)$ is $\ell$-distance-balanced. In Table~\ref{tab:GP} for each pair $(n,k)$ where $5 \leq n \leq 25$ and $2 \leq k < n/2$ the diameter $D$ of $\Gamma = GP(n,k)$ and the set of all $1 \leq \ell \leq D$ for which $\Gamma$ is $\ell$-distance-balanced is given (under the column $\ell$-dist. bal.). The possibility of $k = 1$ is omitted in view of Corollary~\ref{cor:prism} (which is also why we start with $n = 5$). We remark that even though $GP(n,k) \cong GP(n,k^{-1})$ (or $GP(n,-k^{-1})$ if $k^{-1} > n/2$) when $k$ is coprime to $n$ we put both possibilities in the table since one might want to search for patterns just based on the value of $k$.

\begin{table}
\begin{footnotesize}
$$
\begin{array}{|c|c|c||c|c|c||c|c|c|}
\hline
(n,k) & D & \ell-\text{dist. bal.} & (n,k) & D & \ell-\text{dist. bal.} & (n,k) & D & \text{par} \\
\hline
(5,2) & 2 & \{ 1, 2 \} & (6,2) & 4 & \{ 4 \} & (7,2) & 3 & \{ 1, 2, 3 \} \\
(7,3) & 3 & \{ 1, 2, 3 \} & (8,2) & 4 & \{ 4 \} & (8,3) & 4 & \{ 1, 2, 3, 4 \} \\
(9,2) & 4 & \{ 3, 4 \} & (9,3) & 4 & \{ 4 \} & (9,4) & 4 & \{ 3, 4 \} \\
(10,2) & 5 & \{ 1, 2, 3, 4, 5 \} & (10,3) & 5 & \{ 1, 2, 3, 4, 5 \} & (10,4) & 4 & \{ 4 \} \\
(11,2) & 5 & \{ 4, 5 \} & (11,3) & 4 & \{ 4 \} & (11,4) & 4 & \{ 4 \} \\
(11,5) & 5 & \{ 4, 5 \} & (12,2) & 5 & \{ 5 \} & (12,3) & 5 & \{ 5 \} \\
(12,4) & 5 & \{ 5 \} & (12,5) & 4 & \{ 1, 2, 3, 4 \} & (13,2) & 5 & \{ 5 \} \\
(13,3) & 5 & \{ 4, 5 \} & (13,4) & 5 & \{ 4, 5 \} & (13,5) & 4 & \{ 1, 2, 3, 4 \} \\
(13,6) & 5 & \{ 5 \} & (14,2) & 6 & \{ 6 \} & (14,3) & 5 & \{ 5 \} \\
(14,4) & 5 & \{ 1, 2, 3, 4, 5 \} & (14,5) & 5 & \{ 5 \} & (14,6) & 5 & \{ 5 \} \\
(15,2) & 6 & \{ 6 \} & (15,3) & 5 & \{ 1, 4, 5 \} & (15,4) & 5 & \{ 1, 2, 3, 4, 5 \} \\
(15,5) & 5 & \{ 5 \} & (15,6) & 5 & \{ 5 \} & (15,7) & 6 & \{ 6 \} \\
(16,2) & 6 & \{ 6 \} & (16,3) & 6 & \{ 5, 6 \} & (16,4) & 5 & \{ 5 \} \\
(16,5) & 6 & \{ 5, 6 \} & (16,6) & 5 & \{ 3, 4, 5 \} & (16,7) & 5 & \{ 1, 2, 3, 5 \} \\
(17,2) & 6 & \{ 6 \} & (17,3) & 5 & \{ 5 \} & (17,4) & 5 & \{ 1, 2, 3, 4, 5 \} \\
(17,5) & 5 & \{ 1, 2, 3, 4, 5 \} & (17,6) & 5 & \{ 5 \} & (17,7) & 5 & \{ 1, 2, 3, 4, 5 \} \\
(17,8) & 6 & \{ 6 \} & (18,2) & 7 & \{ 7 \} & (18,3) & 6 & \{ 6 \} \\
(18,4) & 5 & \{ 5 \} & (18,5) & 5 & \{ 1, 2, 3, 4, 5 \} & (18,6) & 6 & \{ 6 \} \\
(18,7) & 5 & \{ 1, 2, 3, 4, 5 \} & (18,8) & 6 & \{ 5, 6 \} & (19,2) & 7 & \{ 7 \} \\
(19,3) & 6 & \{ 6 \} & (19,4) & 5 & \{ 5 \} & (19,5) & 5 & \{ 5 \} \\
(19,6) & 6 & \{ 6 \} & (19,7) & 5 & \{ 4, 5 \} & (19,8) & 5 & \{ 4, 5 \} \\
(19,9) & 7 & \{ 7 \} & (20,2) & 7 & \{ 7 \} & (20,3) & 6 & \{ 6 \} \\
(20,4) & 6 & \{ 6 \} & (20,5) & 6 & \{ 6 \} & (20,6) & 6 & \{ 4, 5, 6 \} \\
(20,7) & 6 & \{ 6 \} & (20,8) & 5 & \{ 5 \} & (20,9) & 6 & \{ 1, 2, 3, 5, 6 \} \\
(21,2) & 7 & \{ 7 \} & (21,3) & 6 & \{ 6 \} & (21,4) & 6 & \{ 2, 5, 6 \} \\
(21,5) & 6 & \{ 2, 5, 6 \} & (21,6) & 5 & \{ 5 \} & (21,7) & 6 & \{ 6 \} \\
(21,8) & 5 & \{ 1, 2, 5 \} & (21,9) & 6 & \{ 2, 5, 6 \} & (21,10) & 7 & \{ 7 \} \\
(22,2) & 8 & \{ 8 \} & (22,3) & 7 & \{ 7 \} & (22,4) & 6 & \{ 5, 6 \} \\
(22,5) & 6 & \{ 1, 2, 3, 4, 5, 6 \} & (22,6) & 5 & \{ 5 \} & (22,7) & 7 & \{ 7 \} \\
(22,8) & 6 & \{ 5, 6 \} & (22,9) & 6 & \{ 1, 2, 3, 4, 5, 6 \} & (22,10) & 7 & \{ 5, 6, 7 \} \\
(23,2) & 8 & \{ 8 \} & (23,3) & 6 & \{ 6 \} & (23,4) & 6 & \{ 6 \} \\
(23,5) & 5 & \{ 5 \} & (23,6) & 6 & \{ 6 \} & (23,7) & 6 & \{ 5, 6 \} \\
(23,8) & 6 & \{ 6 \} & (23,9) & 5 & \{ 5 \} & (23,10) & 6 & \{ 5, 6 \} \\
(23,11) & 8 & \{ 8 \} & (24,2) & 8 & \{ 8 \} & (24,3) & 7 & \{ 7 \} \\
(24,4) & 6 & \{ 1, 6 \} & (24,5) & 6 & \{ 1, 2, 3, 4, 5, 6 \} & (24,6) & 6 & \{ 6 \} \\
(24,7) & 6 & \{ 1, 2, 3, 4, 5, 6 \} & (24,8) & 7 & \{ 7 \} & (24,9) & 6 & \{ 5, 6 \} \\
(24,10) & 6 & \{ 5, 6 \} & (24,11) & 7 & \{ 1, 2, 3, 5, 6, 7 \} & (25,2) & 8 & \{ 8 \} \\
(25,3) & 7 & \{ 7 \} & (25,4) & 6 & \{ 6 \} & (25,5) & 6 & \{ 6 \} \\
(25,6) & 6 & \{ 6 \} & (25,7) & 5 & \{ 1, 2, 3, 4, 5 \} & (25,8) & 7 & \{ 7 \} \\
(25,9) & 6 & \{ 6 \} & (25,10) & 6 & \{ 6 \} & (25,11) & 6 & \{ 6 \} \\
(25,12) & 8 & \{ 8 \} & & & & & &\\
\hline
\end{array}
$$
\caption{The $\ell$-distance-balanced property of generalized Petersen graphs.}
\label{tab:GP}
\end{footnotesize}
\end{table}

One of the first things to notice is that each $GP(n,k)$ seems to be $D$-distance-balanced where $D$ is the diameter of $GP(n,k)$. In view of Corollary~\ref{cor:GP} it would suffice to prove that in $GP(n,k)$, where $k \geq 2$, a pair of vertices at diametral distance is always of the form $u_i$, $u_j$ or $v_i$, $v_j$. Unfortunately, this is not the case in general. For instance the graph $GP(7,2)$ is of diameter $3$ but $d(u_0, v_3) = 3$. Nevertheless it does seem that there are not too many pairs $(n,k)$ such that in $GP(n,k)$ there exists some $v_j$ at diametral distance from $u_0$. In fact, a computer search suggests the following might be true. 

\begin{conjecture}
\label{conj:first}
Let $n \geq 3$ and $2 \leq k < n/2$ be integers. If there exists $j \in \ZZ_n$ such that $d(u_0, v_j) = D$, where $D$ is the diameter of $GP(n,k)$, then either $n = 4m$ and $k = 2m-1$ for some $m \geq 3$ or the pair $(n,k)$ is one of $(5,2)$, $(7,2)$ and $(7,3)$.
\end{conjecture}

The reason why we were not able to prove this conjecture in general (we prove that it holds for $k = 2$ in the proof of Theorem~\ref{the:GP}) might be that the diameter and consequently the vertices at diametral distance in $GP(n,k)$ heavily depend on the value of $k$. In fact, to the best of our knowledge, the diameter of the graphs $GP(n,k)$ is not known in general. However, if Conjecture~\ref{conj:first} does hold, then the fenomenon observed in Table~\ref{tab:GP} regarding the diameter does hold in general.

\begin{proposition}
\label{pro:GP}
Suppose that Conjecture~\ref{conj:first} holds, let $n \geq 3$ and $1 \leq k < n/2$ be integers, and let $D$ be the diameter of the generalized Petersen graph $\Gamma = GP(n,k)$. Then $\Gamma$ is $D$-distance-balanced.
\end{proposition}
\proof
By Corollary~\ref{cor:prism} we can assume $k > 1$. Moreover, by Corollary~\ref{cor:GP} we can assume there exists $j \in \ZZ_n$ such that $d(u_0, v_j) = D$, and by assumption that Conjecture~\ref{conj:first} holds, we have that $n = 4m$ and $k = 2m-1$ for some $m \geq 3$ or the pair $(n,k)$ is one of the pairs $(5,2)$, $(7,2)$, $(7,3)$. It is straightforward to check that the graphs $GP(5,2)$ and $GP(7,2) \cong GP(7,3)$ are in fact highly distance-balanced (see also Table~\ref{tab:GP}). 

For the rest of the proof we will thus assume that $n = 4m$ and $k = 2m-1$ for some $m \geq 3$. Since $2k \equiv -2 \pmod{n}$ it is easy to see that for any $0 \leq j \leq 2m$ we have 
$$
	d(u_0, v_j) = \left\{\begin{array}{ccc}
		j+1 & ; & 0 \leq j \leq m, \\
		2m-j+1 & ; & m \leq j < 2m, \\
		3 & ; & j = 2m.\end{array}\right.
$$
Therefore, if some vertex $v_j$ exists, such that $d(u_0, v_j) = D$, it must be that $v_j = v_m$ (or $v_{-m}$). Since $(2m-1)^2 \equiv 1 \pmod {4m}$ it is clear that the permutation $\sigma$ of $V(\Gamma)$, mapping each $u_i$ to $v_{(2m-1)i}$ and each $v_i$ to $u_{(2m-1)i}$ is an automorphism of $\Gamma$. Clearly $\sigma(v_m) = u_{2m^2 - m}$ which is either $u_{-m}$ or $u_m$, depending on whether $m$ is even or odd, respectively. Thus either $\rho^m\sigma$ or $\rho^m\tau\sigma$ interchanges $u_0$ and $v_m$, and so this pair of vertices is balanced by Lemma~\ref{le:swap}. It follows that $\Gamma$ is $D$-distance-balanced.
\qed

The data from Table~\ref{tab:GP} can easily be extended up to at least $n = 200$. The results seem to indicate that for a fixed $k$ there exists some (smallest) integer $n_k$ such that for all $n > n_k$ the graph $GP(n,k)$ is $D$-distance-balanced but is not $\ell$-distance-balanced for any $1 \leq \ell < D$, where $D$ is the diameter of $GP(n,k)$. For instance, it seems that $n_2 = 11$, $n_3 = 16$, $n_4 = 24$, $n_5 = 36$, $n_6 = 48$, $n_7 = 64$, $n_8 = 80$, $n_9 = 100$, $n_{10} = 120$, etc. We therefore make the following conjecture.

\begin{conjecture}
\label{conj:uper}
Let $k \geq 2$ be an integer and let 
$$
	n_k = \left\{\begin{array}{ccc}
		11 & ; & k = 2, \\
		(k+1)^2 & ; & k\ \mathrm{odd},\\
		k(k+2) & ; & k \geq 4\ \mathrm{even}.\end{array}\right.
$$
Then for any $n > n_k$ the graph $GP(n,k)$ is not $\ell$-distance-balanced for any $1 \leq \ell < D$, where $D$ is the diameter of $GP(n,k)$. Moreover, $n_k$ is the smallest integer with this property.
\end{conjecture}

We remark that a result about $1$-distance-balancedness of the graphs $GP(n,k)$, related to Conjecture~\ref{conj:uper}, was proved in~\cite{YHLZ}. In particular, it was proved that for any integer $k \ge 2$ and $n > 6k^2$ the graph $GP(n, k)$ is not 1-distance-balanced (see \cite[Theorem 2]{YHLZ}). 
\medskip

In the reminder of this section we prove that Conjecture~\ref{conj:uper} does hold at least for $k = 2$. We first determine all $1$-distance-balanced $GP(n,2)$ graphs, then all $2$-distance-balanced ones and finally all $\ell$-distance-balanced ones for $\ell \geq 3$. 

For the rest of this section let $\Gamma = GP(n,2)$ for some $n \geq 5$. We first make the following observations regarding the distances in $\G$. Let $0 \leq i \leq n/2$ and consider a shortest path between $u_0$ and $v_i$. Clearly such a path contains just one spoke and at most one outer edge. Observe also that $d(u_0, v_i) = d(v_0, u_i)$ (using the automorphisms $\rho$ and $\tau$). It is thus clear that 
$$
d(u_0, v_i) = d(v_0, u_i) = \left\{\begin{array}{ccc} 1 + \frac{i}{2} & ; & i\ \mathrm{even},\\
								2 + \frac{i-1}{2} & ; & i\ \mathrm{odd}.\end{array}\right.
$$
This enables us to easily calculate the distances between any pair of vertices of $\Gamma$. For instance, if for some $0 \leq i \leq n/2$ every shortest path from $u_0$ to $u_i$ uses at least one inner edge (which clearly occurs if and only if $n \geq 12$ and $i \geq 6$) then we can assume the first edge of such a path is $u_0v_0$, and so $d(u_0,u_i) = d(v_0,u_i)-1$. We also point out that when $n$ is even, every shortest path from $v_0$ to $v_i$ with $i$ even uses only inner edges and is thus of length $i/2$, while every shortest path from $v_0$ to $v_i$ with $i$ odd uses one outer edge and is thus of length $3 + (i-1)/2$. In the case that $n$ is odd, the situation is somewhat different. Namely, in this case, even though one of $(n-1)/2$ and $(n-3)/2$ is odd the shortest path from $v_0$ to the corresponding $v_i$ uses only inner edges (with the first edge being $v_0v_{-2}$), and so this $v_i$ is not closer to $u_0$ than to $v_0$. All this enables us to determine all $1$-distance-balanced generalized Petersen graphs of the form $GP(n,2)$.

\begin{proposition}
\label{pro:GP2_1}
Let $n \geq 5$ be an integer. Then the generalized Petersen graph $GP(n,2)$ is $1$-distance-balanced if and only if $n \in \{5,7,10\}$. 
\end{proposition}
\proof
In view of the automorphisms $\rho$ and $\tau$, Lemma~\ref{le:swap} implies that the graph $\G = GP(n,2)$ is $1$-distance-balanced if and only if the pair $u_0$, $v_0$ is balanced. Using the remarks on distances in $\G$ one can easily determine the sets $W_{u_0 v_0}$ and $W_{v_0 u_0}$ and thus complete the proof. For instance, if $n$ is even then clearly all of the vertices $v_i$ with $i$ even are in $W_{v_0 u_0}$ while all of the vertices $v_i$ with $i$ odd are in $W_{u_0 v_0}$, and so precisely half of the vertices $v_i$ are in $W_{u_0 v_0}$ while the other half is in $W_{v_0 u_0}$. By the above remarks $u_i \in W_{u_0 v_0}$ if and only if $i \in \{0,1,-1\}$ while all the vertices $u_i$ and $u_{-i}$ for $4 \leq i \leq n/2$ are in $W_{v_0 u_0}$. Thus the pair $u_0$, $v_0$ is balanced if and only if $3 = n-7$ that is $n = 10$.  The case when $n$ is odd requires a bit more work but can also be done in a similar way. 

One can first easily check the graphs $GP(5,2)$, $GP(7,2)$ and $GP(9,2)$ by hand (see also Table~\ref{tab:GP}) to verify that out of the three precisely $GP(5,2)$ and $GP(7,2)$ are $1$-distance-balanced. To complete the proof we thus only need to show that if $n \geq 11$ is odd the pair $u_0$, $v_0$ is not balanced. 
As was already pointed out we have $u_i \in W_{u_0 v_0}$ if and only if $i \in \{0,1,-1\}$ and since $u_2$ and $u_3$ are clearly both at equal distances from $u_0$ and $v_0$ we thus also get $u_i \in W_{v_0 u_0}$ if and only if $4 \leq i \leq n-4$. Moreover, $W_{v_0 u_0}$ contains at least all of the vertices $v_i$ and $v_{-i}$ for $0 \leq i \leq n/2$ even, while $W_{u_0 v_0}$ contains the vertices $v_i$ and $v_{-i}$ only for $1 \leq i < (n-3)/2$ odd. It thus follows that 
$$
	|W_{u_0 v_0}| \leq 3 + \frac{n-3}{2} = \frac{n+3}{2}\quad \mathrm{and}\quad  |W_{v_0 u_0}| \geq n-7 + \frac{n-1}{2} = \frac{3n - 15}{2}.
$$
Since in the case that $n \geq 11$ we have $3n-15 > n+3$ this finally shows that for $n \geq 11$ the pair $u_0$, $v_0$ is not balanced. \qed

As it turns out the $1$-distance-balanced graphs $GP(n,2)$ coincide with the $2$-distance-balanced ones.

\begin{proposition}
\label{pro:GP2_2}
Let $n \geq 5$ be an integer. Then the generalized Petersen graph $GP(n,2)$ is $2$-distance-balanced if and only if $n \in \{5,7,10\}$.
\end{proposition}
\proof
In view of the automorphisms $\rho$ and $\tau$ Lemma~\ref{le:swap} implies that $\G = GP(n,2)$ is $2$-distance-balanced if and only if the pairs $u_0$, $v_1$ and $u_0$, $v_2$ are both balanced. Again, one can easily check that for $5 \leq n \leq 14$ both pairs are balanced if and only if $n \in \{5,7,10\}$ (see also Table~\ref{tab:GP}). For the rest of the proof we thus assume $n \geq 14$.

We show that in this case the pair $u_0$, $v_2$ is not balanced. Let $i \in \ZZ_n$ be such that a shortest path from $u_0$ to $u_i$ contains at least one inner edge. Then there also exists a shortest path from $u_0$ to $u_i$ whose second vertex is $v_0$. But since $v_0$ is a neighbor of both $u_0$ and $v_2$, the vertex $u_i$ cannot be closer to $u_0$ than to $v_2$. It is thus clear that $u_i \in W_{u_0 v_2}$ if and only if $i \in \{1,0,-1,-2,-3\}$ (recall that $n \geq 14$). 

Suppose $i \in \ZZ_n$ is such that $v_i \in W_{u_0 v_2}$. Then clearly either $i = 1$ or $i = n - 1 - 2j$ for some small enough $j \geq 0$. If $n$ is even, then the path $(v_2, u_2, u_3, v_3, v_5, v_7, \ldots , v_{n-1-2j})$ is of length $3 + (n-1-2j-3)/2$, and so $2+j < (n+2-2j)/2$ must hold, that is $j < (n-2)/4$. If however $n$ is odd, then the path $(v_2, v_4, \ldots, v_{n-1-2j})$ is of length $(n-1-2j-2)/2$, and so $2+j < (n-3-2j)$ must hold, that is $j < (n-7)/4$. In any case we thus find that (by $\lfloor a \rfloor$ we denote the largest integer not exceeding $a$)
$$
	|W_{u_0 v_0}| \leq 5 + 1 + \left\lfloor \frac{n-3}{4} \right\rfloor + 1 = 6 + \left\lfloor \frac{n+1}{4} \right\rfloor.
$$

Similarly we easily see that $W_{v_2 u_0}$ for sure contains all vertices of the form $v_{2i}$ where $i \geq 1$ and $i-1 < (n-2i)/2+1$ in case $n$ is even and where $i-1< (n-1-2i)/2 + 2$ in case $n$ is odd. In any case $v_{2i} \in W_{v_2 u_0}$ for at least all $1 \leq i \leq \left\lfloor \frac{n+1}{4}\right\rfloor$. It is also not difficult to see that $u_j \in W_{v_2 u_0}$ for at least all $2 \leq j \leq \frac{n+1}{2}$ (for instance, if $n$ is odd then $u_{2i+1} \in W_{v_2 u_0}$ for $i \geq 1$ whenever $i-1+2 < (n-2i-1)/2 + 2$). Thus 
$$
	|W_{v_2 u_0}| \geq \left\lfloor \frac{n+1}{4}\right\rfloor + \frac{n-1}{2}.
$$
For $n \geq 14$ we get $\frac{n-1}{2} > 6$, and so the pair $u_0$, $v_2$ is not balanced.\qed

We are now ready to completely settle the question of $\ell$-distance-balancedness for the graphs $GP(n,2)$. As a consequence we confirm Conjecture~\ref{conj:uper} for $k = 2$.

\begin{theorem}
\label{the:GP}
Let $n \geq 5$, let $\G = GP(n,2)$ and let $D$ be the diameter of $\G$. Then the following holds.
\begin{itemize}
\itemsep = 0pt
\item[(i)] $\G$ is highly distance-balanced if and only if $n \in \{5,7,10\}$.
\item[(ii)] $\G$ is $D$- and $(D-1)$-distance-balanced but not $\ell$-distance-balanced for any $1 \leq \ell \leq D-2$ if and only if $n \in \{9, 11\}$.
\item[(iii)] $\G$ is $D$-distance-balanced but not $\ell$-distance-balanced for any $1 \leq \ell \leq D-1$ if and only if $n \notin \{5,7,9,10,11\}$. 
\end{itemize} 
In particular, if $n > 11$, then $\G$ is $D$-distance-balanced but is not $\ell$-distance-balanced for any $1 \leq \ell \leq D-1$.
\end{theorem}
\proof
As in the previous two proofs the argument is much easier for large values of $n$, so we verify the cases with small $n$ separately. We can thus verify that the statement of the theorem is true for all $n  \leq 12$ (see also Table~\ref{tab:GP}). For the rest of the proof we thus assume $n \geq 13$. Observe that this implies $d(u_0, u_6) = 5$, and so $D \geq 5$.

We first prove that $\G$ is $D$-distance-balanced. We establish this by proving that there exists no $v_i$ such that $d(u_0, v_i) = D$. That $\G$ is $D$-distance-balanced then follows from Corollary~\ref{cor:GP}. Suppose to the contrary that such a vertex $v_i$ exists and assume with no loss of generality that $i \leq n/2$. Since $D$ is the diameter of $\Gamma$ we have $d(u_0, u_i) \in \{D, D-1\}$. Now, $d(u_0, v_5) = d(u_0, v_6) = 4 < 5 = d(u_0, u_6)$, and so $i > 6$. Consequently, every shortest path $P$ from $u_0$ to $u_i$ uses at least one inner edge. But then there clearly must also exist a shortest path from $u_0$ to $u_i$, whose last edge is $v_i u_i$. However, this implies that $d(u_0, v_i) < d(u_0, u_i)$, a contradiction, which thus proves that $\G$ is $D$-distance-balanced. 

To complete the proof we now only need to show that for any $1 \leq \ell < D$ the graph $\G$ is not $\ell$-distance-balanced. Propositions~\ref{pro:GP2_1} and~\ref{pro:GP2_2} show that this is true for $\ell \in \{1,2\}$, and so we can assume $3 \leq \ell \leq D-1$. We first show that there exists $v_i$ such that $d(u_0, v_i) = \ell$. Indeed, by the above argument no $v_j$ is at distance $D$ from $u_0$ so for some $j \leq n/2$ we have that $d(u_0, u_j) = D$ or $d(v_0, v_j) = D$. In any case $d(u_0, v_j) = D-1$. But a shortest path from $u_0$ to $v_j$ for sure uses only one spoke and at most one outer edge, and so there is a shortest path $P$ from $u_0$ to $v_j$ of length $D-1$ such that except for perhaps the first two vertices all of its vertices are of  the form $v_i$. But then the vertex on $P$ preceeding $v_j$ is of the form $v_i$ and is at distance $D-2$ from $u_0$. The one before it is also in $\{v_i \colon i \in \ZZ_n\}$ and is at distance $D-3$ from $u_0$, etc. 

\begin{figure}[!h]
\begin{center}
\includegraphics[scale=0.8]{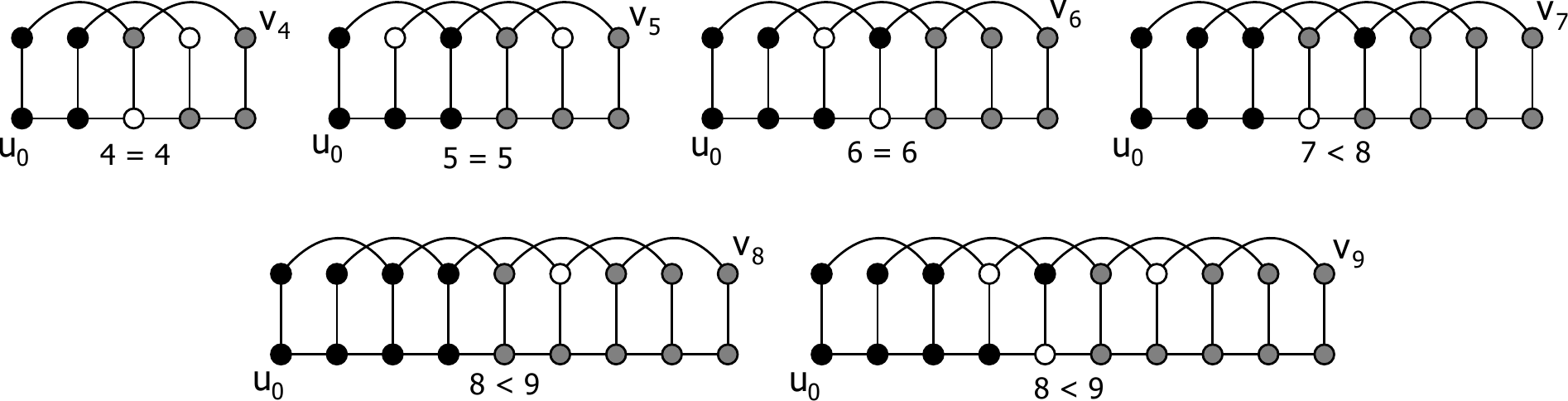}
\end{center}
\caption{The graphs $\G_1$ for $4 \leq i \leq 9$.}
\label{fig:strip2}
\end{figure}

\begin{figure}[!h]
\begin{center}
\includegraphics[scale=0.8]{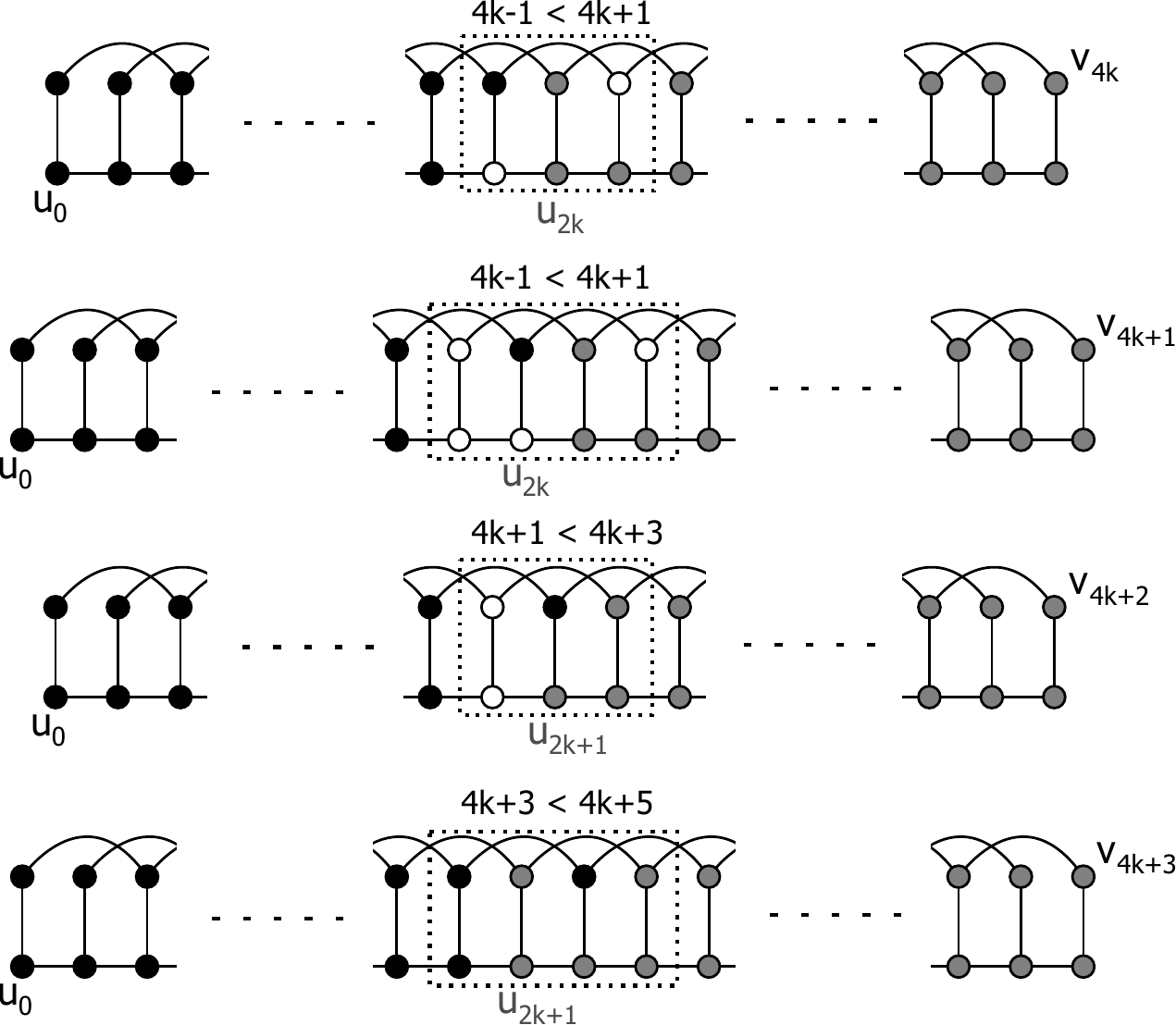}
\end{center}
\caption{The graphs $\G_1$ for $i \geq 10$.}
\label{fig:strip3}
\end{figure}

Let now $v_i$ with $i \leq n/2$ be such that $d(u_0, v_i) = \ell$. Since $\ell \geq 3$ and $d(u_0, v_3) = d(u_0, v_4) = 3$ we can assume that $4 \leq i \leq n/2$. Let now $V_1 = \{u_j \colon 1 \leq j \leq i-1\} \cup \{v_j \colon 1 \leq j \leq i-1\}$ and $V_2 = \{u_j \colon i+1 \leq j \leq n-1\} \cup \{v_j \colon i+1 \leq j \leq n-1\}$. It is now clear that for any $0 \leq j \leq i$ every shortest path from either $u_0$ or $v_i$ to either $u_j$ or $v_j$ uses no vertex from $V_2$ and similarly for any $i \leq j \leq n$ every shortest path from either $u_0$ or $v_i$ to either $u_j$ or $v_j$ uses no vertex from $V_1$. To determine the sets $W_{u_0 v_i}$ and $W_{v_i u_0}$ we can thus separately consider the graph $\G_1 = \G - V_2$, obtained from $\G$ by deleting all the vertices from $V_2$, and the graph $\G_2 = \G - V_1$. The graphs $\G_1$ and $\G_2$ both have a similar structure. They only differ in their order (which for both is at least $10$). It thus suffices to analyze all of the possibilities for $\G_1$. The situation regarding which vertices of $\G_1$ are closer to $u_0$ than to $v_i$ and vice versa is somewhat different for $i \leq 9$ than for $i \geq 10$ when it only depends on the congruence of $i$ modulo $4$. We present all of the possibilities for $4 \leq i \leq 9$ on Figure~\ref{fig:strip2} and the possibilities for $i \geq 10$, depending on the congruence of $i$ modulo $4$, on Figure~\ref{fig:strip3}. 
We find that in $\G_1$ we always have $|W_{v_i u_0}| \geq |W_{u_0 v_i}|$ and moreover, equality holds only for $4 \leq i \leq 6$. Since a similar situation holds for $\G_2$ we see that if either $i \geq 7$ or $n-i \geq 7$ the pair $u_0$, $v_i$ is not balanced in $\G$. But if $i \leq 6$ and $n-i \leq 6$, then $n \leq 12$, a contradiction. Thus the pair $u_0$, $v_i$ is not balanced in $\G$, and so $\G$ is indeed not $\ell$-distance-balanced, as claimed. This completes the proof.
\qed

\section{Suggestions for further research}
\label{sec:FR}

We conclude the paper with some suggestions for future research. In~\cite{Ha} Handa proved that every $1$-distance-balanced graph is $2$-connected. In his PhD thesis~\cite{F} Frelih proved that this is no longer the case if we move to $2$-distance-balanced graphs and characterized connected $2$-distance-balanced graphs, which are not $2$-connected. We therefore propose the following problem.

\begin{problem}
For each $\ell \geq 3$ characterize all connected $\ell$-distance-balanced graphs which are not $2$-connec\-ted.
\end{problem}

In~\cite{Ha} Handa then asked whether all bipartite $1$-distance-balanced graphs were also $3$-connected. A negative answer to this question was given in~\cite{MS} with an infinite family of examples. It turned out that even though bipartite $1$-distance-balanced graphs which are not $3$-connected exist, they have a rather restricted structure. The next problem is thus the following.

\begin{problem}
Generalize the results of~\cite{MS} to the class of bipartite $\ell$-distance-balanced graphs, $\ell \ge 2$, which are not $3$-connected.
\end{problem}

Recently two subfamilies of $1$-distance-balanced graphs were introduced, namely the strongly distance-balanced graphs~\cite{KMMM06} (which coincide with the distance degree regular graphs) and the nicely distance-balanced graphs~\cite{KM}. These two concepts could easily be extended to $\ell$-distance-balanced graphs to obtain strongly $\ell$-distance-balanced graphs and nicely $\ell$-distance-balanced graphs. We thus propose the following problem.

\begin{problem}
Introduce the concepts of strongly and nicely $\ell$-distance-balanced graphs and investigate the properties of such graphs. 
\end{problem} 

In the past few years various results describing how the $1$-distance-balanced property (and strongly $1$-distance-balanced property) of graphs is preserved under various graph products (see for instance~\cite{BCPSSS, JKR}). Moreover, in~\cite{F} Frelih investigated $2$-distance-balanced graphs with respect to Cartesian and lexicographic products. We thus propose to study these things more generally. 

\begin{problem}
Study $\ell$-distance-balanced graphs with respect to various graph products.
\end{problem}

In Section~\ref{sec:D<4} the $\ell$-distance-balancedness for graphs of diameter at most $3$ was investigated where for diameter $3$ we restricted ourselves to bipartite graphs. We believe some interesting results could be obtained also for nonbipartite graphs of diameter $3$, as well as for bipartite graphs of diameter $4$.

\begin{problem}
Generalize the results of Section~\ref{sec:D<4} to non-bipartite graphs of diameter 3 and to (bipartite) graphs of diameter $4$.
\end{problem}

In this paper we also considered the $\ell$-distance-balanced property for cubic graphs in some detail. However, as we saw in Section~\ref{sec:GP} even the generalized Petersen graphs, which are a very special subfamily of cubic graphs, seem to present quite a hard problem when it comes to $\ell$-distance-balancedness. The two conjectures from Section~\ref{sec:GP} (and the related Proposition~\ref{pro:GP}) are just two problems regarding the $GP(n,k)$ graphs and their $\ell$-distance-balancedness that can be considered. There is at least one other interesting problem regarding these graphs that should be considered. Upon inspection of Table~\ref{tab:GP} one quickly notices that there are not so many pairs $(n,k)$ for which the graph $GP(n,k)$ is highly distance-balanced but it seems there are infinitely many such pairs. It is thus very natural to consider the following problem. 

\begin{problem}
Determine all pairs of integers $(n,k)$, where $n \geq 5$ and $2 \leq k < n/2$, such that the generalized Petersen graph $GP(n,k)$ is highly distance-balanced or at least determine whether there are infinitely many such pairs.
\end{problem}

Recall that we proved in Proposition~\ref{pro:abelian} that every connected Cayley graph of an abelian group is highly distance-balanced. As was pointed out in Section~\ref{sec:basic} this is not true for all Cayley graphs (see the example from Figure~\ref{fig:examples}). However, one might get some similar results for Cayley graphs over groups which are ``close'' to being abelian, say dihedral groups. We thus propose to study the $\ell$-distance-balancedness property of Cayley graphs of dihedral groups.

\begin{problem}
For each connected Cayley graph $\G$ of a dihedral group determine all $\ell \geq 1$ such that $\G$ is $\ell$-distance-balanced. If this is too difficult in general, consider this problem at least for cubic Cayley graphs of dihedral groups.
\end{problem}


\end{document}